\documentclass{amsproc}
\usepackage{amssymb}
\usepackage{amsmath}
\usepackage{amsfonts}

\newtheorem{theorem}{Theorem}[section]

\newtheorem{definition}[theorem]{Definition}

\newtheorem{proposition}[theorem]{Proposition}

\theoremstyle{remark}

\numberwithin{equation}{section}

\newcommand{\R}{\mathbb{R}}

\newcommand{\func}[1]{\mathop{\mbox{\rm #1}}}

\begin{document}

\title{EQUIVALENCE\ FOR\ DIFFERENTIAL\ EQUATIONS}
\author{Odinette Ren\'{e}e Abib}
\address{Laboratoire de Math\'{e}matiques Rapha\"{e}l Salem\\
UMR 6085 CNRS, Universit\'{e} de Rouen\\
76801, Saint Etienne du Rouvray}
\email{renee.abib@univ-rouen.fr}



\begin{abstract}
We shall study the equivalence problem for ordinary differential
equations with respect to the affine transformation group 
$A( 2,\R)$.
\end{abstract}

\maketitle

\section{Introduction}

Two differential equations are called equivalent if one can be transformed
into another by a certain change of variables. In particular this change of
variables transforms solutions of one equation into solutions of another. In
general, the problem of equivalence of differential equations consists in
determining whether two equations are equivalent up to a given class of
transformations. W. Kry\'{n}sky \cite{Kr} and B. Dubrov \cite{DKM} examine
differential equations up to contact transformations. Sophus Lie was the
first to use some approach to the problem of equivalence of differential
equations (description of invariants, computation of symetry group).

In \cite{Fe} Fels consider the problem of equivalence between two systems of
second-order differential equations

\begin{eqnarray}
\frac{d^{2}x^{i}}{dt^{2}} &=&f^{i}(t,x^{j},\frac{dx^{j}}{dt})\text{ \ \ \ \ }%
(1\leq i,\text{ }j\leq n)\text{\ }\text{\ } \label{eq:fle}\\
\frac{d^{2}y^{i}}{ds^{2}}\text{\ } &=&g^{i}(s,y^{j},\frac{dy^{j}}{ds})\text{
\ \ \ \ }(1\leq i,\text{ }j\leq n)\text{ \ \ }\label{eq:chern} \\
&&\text{\ \ }\notag
\end{eqnarray}%
under the pseudo-group of smooth invertible local point transformations

\begin{equation*}
\psi \left( t,x^{j}\right) =\left( s,y^{i}\right) =\left( \phi \left(
t,x^{j}\right) ,\varphi _{i}\left( t,x^{j}\right) \right)
\end{equation*}

This notion of two system being equivalent defines an equivalence relation
on the set of differential equations on the form \eqref{eq:fle}. Fels was able to
cast the question of equivalence between \eqref{eq:fle} and \eqref{eq:chern} into a question
about the equivalence of exterior differential systems associated on jet
space $J^{1}\left( 
\mathbb{R}
,%
\mathbb{R}
^{n}\right) ,$ at wich point the Cartan's \ method \cite{Ca1, EKMR, Ga, GS, Ol2, QT}                                             be may applied. The problem of equivalence for $n=1$ was
originally solved by Cartan in \cite{Ca1}. Chern in \cite{Ch} investigated the two
equivalence problems for systems under the restricted pseudo-groups of
smooth invertible local transformations wich preserve the independent
variable as given by

\begin{equation*}
\psi \left( t,x^{j}\right) =\left( t,\varphi _{i}\left( x^{j}\right) \right) 
\text{ and }\psi \left( t,x^{j}\right) =\left( t,\varphi _{i}\left(
t,x^{j}\right) \right) \text{ }
\end{equation*}

In the previous papers \cite{Ab1, Ab2} the relationship between differential
equations, Pfaffian systems and geometric structures are study. We have seen
that every differential equation can be expressed as a Pfaffian system
satisfying the structure equation and that the integration of a given
equation is deeply related to the structure equation. We shall show it by
means of interesting examples. My contribution here is the study of
equivalence problem for the family of ordinary equations with respect to the
affine transformation group $A\left( 2,%
\mathbb{R}
\right) $ (section 4). There exist 19 different types of first order
ordinary differential equations which admit at least 1- dimensional Lie
groups in $A\left( 2,%
\mathbb{R}
\right) $. All the equations which belong to the above types can be
integrated by quadrature.

In the paper, by the word differentiable we mean always differentiable of
classe $C^{\infty }$.

\textbf{Acknowledgments. }Special thanks go to Witold Respondek and Yuri
Sachkov for his encouragement and support.

\section{Basic Definition, Examples}

\subsection{\textbf{Pfaffian systems}}

In this subsection we will review some basic concepts and facts on Pfaffian
systems theory \cite{BCGGG}. Let $M$ be an differentiable manifold. $F(M)$
denotes the ring of real-valued differentiable functions on $M$ and 
${ %
\Lambda }^{1}(M)$ the $F(M)$-module of all 1-forms (Pfaffian forms) on $M$%
. \ 

A $F(M)$-submodule $\Sigma$ of ${\Lambda }^{1}(M)$ is called a 
\textit{Pfaffian system of rank }$n$\textit{\ on M}\textbf{\ }if $\Sigma $
is generated by $n$ linearly independent Pfaffian forms $\theta ^{1},\ldots
,\theta ^{n}$. A submanifold $N$ of $M$ is said an \textit{integral manifold 
}of $\Sigma $ if $i^{\ast }(\theta )=0$ for all $\theta \in \Sigma $, where $%
i$ denotes the immersion $N$ $\hookrightarrow $ $M$. A differentiable
function is said a \textit{first integral of }$\Sigma $ if the exterior
derivate $df$ belongs to $\Sigma .$ By the symbol $\Sigma =\langle \theta
^{1},\ldots ,\theta ^{n}\rangle $ we mean that the Pfaffian system $\Sigma $
is generated by linearly independent Pfaffian forms $\theta ^{1}\ldots
,\theta ^{n}$ defined on $M$.

For each Pfaffian systems $\Sigma$ on $M$ we can construct the dual sytem,
that is, the differentiable subbundle $D=D(\Sigma )$ of the tangent bundle $%
T(M)$ on $M$ such that the fiber dimension of $D$ is equal to $\dim (M)-n$.
Let $\Gamma (D)$ be the sheaf of germs of local vector fields which belong
to $D$ and $\Gamma (D)_{x},x\in M,$ the stalk of $\Gamma (D)$ at $x$. For $%
x\in M$, we defined the subspaces $Ch(D)_{x\text{ }}$of $T_{x}(M)$ by

\begin{equation*}
Ch(D)_{x}=\left\{ X_{x}\in D_{x};\left[ X_{x}^{\prime },\Gamma (D)_{x}\right]
\subseteq \Gamma (D)_{x}\right\}
\end{equation*}

where $X$ denotes a vector field and $X_{x}^{\prime }$ the germ at $x$
determined by $X$. We suppose that $\dim Ch(D)_{x}$ is constant on $M$.
Thus, we obtain the subbundle\ $Ch(D)$ of $T(M)$ called \ \textit{Cauchy
characteristic} of $D$. The distrubution characteristic of $D$ is the module
spanned by all vector fields $Y$ \ that belongs to $D$ such that $\left[ Y,D%
\right] \subseteq D$. The dual system $Ch(\Sigma )$ of $Ch(D)$ is called the 
\textit{Cauchy characteristic system} of $\Sigma$. The following theorem is
due to E. Cartan. \ 

\begin{theorem}\label{theo:cartan}
Let $\Sigma =\langle \theta ^{1},\ldots ,\theta ^{n}\rangle $ be a Pfaffian
system.

1. If $\Sigma $ is completely integrable, i.e. $d\theta ^{i}=0(\func{mod}%
\theta ^{1},\ldots ,\theta ^{n})$ ,\ $\left( 1\leq i\leq n\right) ,$ then $%
Ch(\Sigma )=\Sigma$. \ \ \ \ \

2. If $\Sigma $ \textit{is not completely integrable, then there exist
linearly independent Pfaffian forms }$w^{1},\ldots ,w^{m}$ \textit{%
satisfying the following conditions:}\ \ 

\ \ \ \ \ \ $(i)$ $\theta ^{1},\ldots ,\theta ^{n},w^{1},\ldots ,w^{m}$ 
\textit{are also linearly independen}t;\ \ 

\ \ \ \ \ \ $(ii)$ $(\theta ^{1},\ldots ,\theta ^{n},w^{1},\ldots ,w^{m})$ 
\textit{forms a (local) generator of }$Ch(\Sigma )$;

\ \ \ \ \ $(iii)$ $d\theta^{i}=
\displaystyle{\sum\limits_{j,k=1}^{m}}C_{jk}^{i}w^{i}\wedge w^{k}$ $\left( \func{mod}%
\theta ^{1},\ldots ,\theta ^{n}\right) ,$ \textit{where }$C_{jk}^{i}$ 
\textit{denotes a differentiable function }$($\textit{1}$\leq i\leq n,$%
\textit{\ 1}$\leq j,k\leq m).$\textit{\ }

\ \textit{3. }$Ch(\Sigma )$ \textit{is completely integrable. \ }

\ \textit{4. Let }$x^{1},\ldots ,x^{n+m}$ \textit{be independent first
integrals of} $Ch(\Sigma ).$ \textit{Then there exist linearly independent
Pfaffian forms} $\eta ^{i}={\displaystyle\sum\limits_{j=1}^{n+m}}A_{j}^{i}(x^{1},\ldots
,x^{n+m})dx^{j},$ $i=1,\ldots ,n,$ \textit{such that} $(\eta ^{1},\ldots
,\eta ^{n})$\ \ forms a (local) generator of $\Sigma$.
\end{theorem}

By making use of the property 2, we can construct the Cauchy characteristic
system $Ch(\Sigma)$.

\begin{definition}
\textit{A system }$(w^{1},\ldots ,w^{m})$ \textit{of linearly independent
Pfaffian forms on }$M$ \textit{will be said a solvable system of }$\Sigma
=\langle \theta ^{1},\ldots ,\theta ^{n}\rangle $ \textit{if it satisfies
the following conditions:}

\textit{(i) }$(w^{1},\ldots ,w^{m})$ \textit{forms a generator of }$%
Ch(\Sigma ).$

(ii) $dw^{1}=0$ \textit{and }$dw^{p}=0\func{mod}(w^{1},\ldots ,w^{p-1})$ for
all $p=2,\ldots ,m.$
\end{definition}

If we can find a solvable system of $\Sigma ,$ then $m$ independents first
integral of $Ch(\Sigma )$ are given par quadrature.

\subsection{Examples}

In the subsection we shall consider by means of examples \cite{Ab2} the relation
between the differential equations, Pfaffian systems and structure equation.

\textbf{a) }Consider the Pfaffian system $\Sigma =\langle \theta \rangle ,$ $%
\theta =dz+pdx+p^{2}dy,$ on $%
\mathbb{R}
^{4}=\left\{ \left( x,y,z,p\right) \right\} $. We have $d\theta =dp\wedge
\left( dx+2pdy\right) $ and 
\begin{equation*}
w^{1}=dp,\text{ }w^{2}=dx+2pdy,\text{ }w^{3}=\theta
\end{equation*}

determine the Cauchy characteristic system of $\Sigma .$ We can find by
quadrature three independent first integrals as follows:

\begin{equation*}
u_{1}=z+xp+yp^{2},\text{ }u_{2}=x+2yp,\text{ }u_{3}=p
\end{equation*}

and $\theta $ itself is expressed as $\theta =du_{1}-u_{2}du_{3}.$ The
system $(w_{1},w_{2},w_{3})$ is a solvable system of $\Sigma .$

\textbf{b) }We consider an absolute parallelism $%
w^{1},w^{2},w^{3},w^{4},w^{5},w^{6}$ on $%
\mathbb{R}
^{6}$ satisfying the equations

\begin{eqnarray}
dw^{1} &=&0,\text{ }dw^{2}=0\func{mod}(w^{1},w^{2}), \notag\\
dw^{3} &=&w^{1}\wedge w^{4}+w^{2}\wedge w^{5}\text{ }\func{mod}(w^{3}),
   \label{eq:parallelism}
    \text{\ \ } \\
dw^{4} &=&0\func{mod}(w^{3},w^{4},w^{5}), \notag\\
dw^{5} &=&w^{2}\wedge w^{6}\func{mod}(w^{3},w^{4},w^{5}).\notag
\end{eqnarray}

Let $x$ and $y$ be two independent first integrales of the completely
Pfaffian integrable system $w^{1}=w^{2}=0;$ the form $w^{3}$ is expressed as

\begin{equation*}
w^{3}=a(dz-pdx-qdy),\text{ \ \ }a\neq 0.
\end{equation*}

The functions $x,y,z,p$ and $q$ are independent first integrals of the
completely integrale Pfaffian system $w^{1}=w^{2}=w^{3}=w^{4}=w^{5}=0.$%
Therefore $w^{4}$ and $w^{5}$ can be written by means of the exterior
derivates $dx,dy,dz,dp,dq$ and the formulas

\begin{eqnarray*}
dp-rdx-sdy &=&a_{1}w^{4}+a_{2}w^{5}+a_{3}w^{3} \\
dq-s^{\prime }dx-tdy &=&a_{4}w^{4}+a_{5}w^{5}+a_{6}w^{3}
\end{eqnarray*}

determine the fonctions $r,s,s^{\prime },t$ and $a_{i}^{\prime }$ $s$ of the
variables $x,y,z,p,q$ and another $u.$ From the equation $dw^{3}=w^{1}\wedge
w^{4}+w^{2}\wedge w^{5}(\func{mod}w^{3}),$ one can verify that the function $%
s$ coincides with $s^{\prime }.$ Moreover, the equations $dw^{4}=0,$ $%
dw^{5}=w^{2}\wedge w^{6}\func{mod}(w^{3},w^{4},w^{5}) $ imply

\begin{equation*}
rank{\Bigl (}\frac{\partial r}{\partial u},\frac{\partial s}{\partial u},%
\frac{\partial t}{\partial u}{\Bigr )=1.}
\end{equation*}

Therefore the functions

\begin{equation*}
r=r(x,y,z,p,q,u)\text{ \ },\text{ \ }s=s(x,y,z,p,q,u)\text{ },\text{ }%
t=t(x,y,z,p,q,u)
\end{equation*}

determine a system of second-order partial differential equations. This
family of systems of differential equations determined by an absolute
parallelism satisfying \eqref{eq:parallelism} is the main subject of Cartan's researches in
his paper \cite{Ca1,Ca2}. For example, take the system of differential
equations

\begin{equation}
 \label{eq:surface}
\frac{\partial ^{2}z}{\partial x^{2}}=0\text{ },\text{ }\frac{\partial ^{2}z%
}{\partial x\partial y}=z-x\frac{\partial z}{\partial x}
\end{equation}

Putting on $%
\mathbb{R}
^{6}=\left\{ (x,y,z,p,q,t)\right\} ,$ $w^{1}=dx,$ $w^{2}=dy,$ $%
w^{3}=dz-pdx-qdy,$ $w^{4}=dp-(z-xp)dy,$ $w^{5}=dq-(z-xp)dx-tdy$ \ and $%
w^{6}=dt-(q-x(z-xp))dx$ we have the structure equations

\begin{eqnarray*}
dw^{1} &=&0,\text{ }dw^{2}=0, \\
dw^{3} &=&w^{1}\wedge w^{4}+w^{2}\wedge w^{5}, \\
dw^{4} &=&w^{2}\wedge w^{3}-xw^{2}\wedge w^{4}, \\
dw^{5} &=&w^{2}\wedge w^{6}+w^{1}\wedge w^{3}-xw^{1}\wedge w^{4}, \\
dw^{6} &=&w^{1}\wedge w^{5}-xw^{1}\wedge w^{3}-x^{2}w^{1}\wedge
w^{4}+Kw^{1}\wedge w^{2},
\end{eqnarray*}

where $K=t-xq+x^{2}(z-xp).$ The absolute parallelism satisfies the equations
\eqref{eq:parallelism}. It easy to see that the system $\left(
w^{2},w^{3},w^{4},w^{5},w^{6}\right) $ forms a solvable system of $\Sigma
=\langle w^{3},w^{4},w^{5}\rangle .$ Five independent first integrals of the
solvable system \ are given by quadrature as follows:

\begin{equation*}
u_{1}=y,\text{ }u_{2}=z-xp,\text{ }u_{3}=p,\text{ }u_{4}=q-x(z-xp),\text{ }%
u_{5}=K\text{ },
\end{equation*}

and we have \cite{Ca2}

\begin{eqnarray*}
w^{3}-xw^{4} &=&du_{2}-u_{4}du_{1}, \\
w^{4} &=&du_{3}-u_{2}du_{1,} \\
w^{5}-xw^{3} &=&du_{4}-u_{5}du_{1}.
\end{eqnarray*}

By the expression, the general integral surface of \eqref{eq:surface} is given by the
formulas:

\begin{equation*}
p=f(y),\text{ }z-xp=f\text{ }^{\prime }(y),\text{ }q-x(z-xp)=f\text{ }%
^{\prime \prime }(y),\text{ }t-x(q-x(z-xp))=f\text{ }^{\prime \prime \prime
}(y)
\end{equation*}

where $f$ is a differentiable function and $f$ $^{\prime },$ $f$ $^{\prime
\prime }$and $f$ $^{\prime \prime \prime }$ denote its derivates.

\section{Equivalence}

In the paper \cite{Ab2} , we have seen that every differential equation can be
expressed as a Pfaffian system satisfying the structure equation and the
integration of a given equation is deeply related to the structure equation.
In this section, we shall consider the equivalence problem for Pfaffian
systems and hence for differential equations under the action of Lie groups.

Let $M$ be a differentiable manifold, $G$ be a Lie group acting on $M$ on
the left. For a Pfaffian system $\Sigma$ on $M$ we set

\begin{equation*}
g^{\ast }\Sigma =\left\{ L_{g}^{\ast }\theta \in \Lambda ^{1}M;\text{ }%
\theta \in \Sigma \right\} ,
\end{equation*}

and%
\begin{equation*}
G(\Sigma )=\left\{ g\in G;\text{ }g^{\ast }\Sigma =\Sigma \text{ }\right\} ,
\end{equation*}

where $L_{g}$ denotes the left action of $g\in G$ on $M$; $g^{\ast }\Sigma $
is a Pfaffian system on $M$ and $G(\Sigma )$ is a subgroup of $G$ .

Two Pfaffian systems $\Sigma _{1}$ and $\Sigma _{2}$ on $M$ are \textit{%
equivalent under the action of }$G$ if there is an element $g\in G$ such
that $g^{\ast }\Sigma _{1}=\Sigma _{2}$.

Let $F$ be a family of Pfaffian systems on $M$. The problems (Lie
programme) to be solved are as follows:

\textbf{1) }Determine the condition for the equivalence of the elements of
$F$.

\textbf{2) }Classify the Pfaffian systems in $F$ under the action
of $G$.

\textbf{3) }For each $\Sigma \in F$, determine the structure of
the subgroup $G(\Sigma )$.

\textbf{4) }Research the relation between the integration of a Pfaffian
system $\Sigma$ and the structure of the subgroup $G(\Sigma )$, $i.e.$,
reduce the integration of a Pfaffian system $\Sigma$ to auxiliary systems
obtained via the knowledge of the structure of $G(\Sigma )$.

\textit{In section 4 we consider the problem 2) in the particular case where 
}$G=A(2,%
\mathbb{R}
).$

For a Lie subgroup $G^{\prime }$ of $G$, we set

\begin{equation}
\label{eq:greno}                                                                                                                 F(G') =\left\{ \Sigma \in F ;\text{ }G(\Sigma
)=G^{\prime }\right\}. 
\end{equation}

For every $\Sigma \in F(G')$, $G'$ is the
largest subgroup of $G$ which leaves $\Sigma $ invariant. It is possible
that $F(G')$ is an empty set.

\begin{proposition}
Let $\Sigma ,$ $\Sigma _{1},$ $\Sigma _{2}$ be Pffafian systems $\in F$
and let $G^{\prime }$ be a subgroup of $G$.

1. For any $g\in G,$ $G(g^{\ast }\Sigma )=g^{-1}G(\Sigma )g$.

2. If $gG(\Sigma _{2})g^{-1}=G(\Sigma _{1})=G^{\prime }$ for an element $%
g\in G,$ then $\Sigma _{1\text{ }}$and $g^{\ast }\Sigma _{2}$ are in
$F(G')$.

3. The normalizer $N(G^{\prime }:G)$ of $G^{\prime }$ in $G$ acts on
$F(G')$.

4. If $\Sigma _{1},$ $\Sigma _{2}\in F(G')$ and $g^{\ast
}\Sigma _{2}=\Sigma _{1}$ for an element $g\in G,$ then $\Sigma _{1}$and $%
\Sigma _{2}$ lie in the same orbit determined by the action of $N(G^{\prime
}:G)$ on $F(G')$.
\end{proposition}

\textit{Proof.} 1) and 2) may be clear. To prove 3), suppose that $\Sigma
\in F(G')$ and $g\in N(G^{\prime }:G)$. From \eqref{eq:greno} and
1) of the proposition, we have

\begin{equation*}
G(g^{\ast }\Sigma )=g^{-1}G(\Sigma )g=g^{-1}G^{\prime }g=G^{\prime }
\end{equation*}

and hence $g^{\ast }\Sigma \in F(G')$.

4) From \eqref{eq:greno} and 1) of proposition, we have%
\begin{equation*}
g^{-1}G^{\prime }g=g^{-1}G(\Sigma _{2})g=G(g^{\ast }\Sigma _{2})=G(\Sigma
_{1})=G^{\prime }
\end{equation*}

and hence $g\in N(G':G)$.\

By virtue of this proposition, the equivalence probleme and the
classification are reduced to the following problems:

\textit{(i) }determine all conjugate classes of the subgroups of $G$.

\textit{(ii) }For a representative\textit{\ }$G^{\prime }$ of each conjugate
class, determine the set $F(G')$.

\textit{(iii) }Observe the action of $N(G^{\prime }:G)$ on $F(G')$. \ 

Since there are, in general, many subgroups $G^{\prime }$ of $G$ such that $%
F(G')$ is empty set, this reduction of the problems is
not always the best one. Moreover, the Pfaffian systems to be considered are
not always defined globally on $M$\ . Therefore, instead of ordinary Lie
groups, we have to consider Lie pseudogroups \cite{Ca1, KR, Ol1, Ol3}.
Then the subject of the study is {\em invariants} of a Pfaffian system with respect to a given
Lie pseudogroup. At the rate, we can recognize that the subgroup $G(\Sigma )$%
\ plays an important role in the problems.\ \ \ \ \ \ \ \ \ \ \ \ \ \ \ \ \

\section{Equivalence with respect to the A(2,$%
\mathbb{R}
$)}

Let $\ G$ be a finite dimensional Lie group and let $\Sigma$ be a
left-invariant completely integrable Pfaffian system on $G$. We denote by $%
I_{g}(\Sigma )$ the maximal integral manifold through $g\in G$ and we set

\begin{equation*}
G_{g}(\Sigma )=\left\{ h\in G;\text{ }L_{h}(I_{g}(\Sigma ))=I_{g}(\Sigma
)\right\} \text{ },
\end{equation*}

Since $\Sigma$ is left-invariant, $G_{g}(\Sigma ),$ $g\in G,$ are mutually
conjugate in $G$.

\subsection{Invariant forms of A(2,$%
\mathbb{R}
$)}

Let $A(2,%
\mathbb{R}
)$ be the affine transformation group on $%
\mathbb{R}
^{2}.$ By making use of the matrix representaion

\begin{equation*}
A(2,%
\mathbb{R}
)=\left\{ 
\begin{bmatrix}
x_{3} & x_{4} & x_{1} \\ 
x_{5} & x_{6} & x_{2} \\ 
0 & 0 & 1%
\end{bmatrix}%
;\text{ }x_{3}x_{6}-x_{4}x_{5}\neq 0,\text{ }x_{i}\in 
\mathbb{R}
,\text{ }i=1,2,\ldots ,6\right\} ,
\end{equation*}

we have a basis of invariant forms of $A(2,%
\mathbb{R}
)$

\begin{eqnarray*}
w^{1} &=&\frac{1}{D}\left( x_{6}dx_{1}-x_{4}dx_{2}\right) , \\
w^{2} &=&\frac{1}{D}\left( x_{3}dx_{2}-x_{5}dx_{1}\right) , \\
w^{3} &=&\frac{1}{D}\left( x_{6}dx_{3}-x_{4}dx_{5}\right) , \\
w^{4} &=&\frac{1}{D}\left( x_{6}dx_{4}-x_{4}dx_{6}\right) , \\
w^{5} &=&\frac{1}{D}\left( x_{3}dx_{5}-x_{5}dx_{3}\right) , \\
w^{6} &=&\frac{1}{D}\left( x_{3}dx_{6}-x_{5}dx_{4}\right) ,
\end{eqnarray*}

where we put $D=x_{3}x_{6}-x_{4}x_{5}.$ We have then the structure equation

\begin{eqnarray}
dw^{1} &=&w^{1}\wedge w^{3}+w^{2}\wedge w^{4}, \notag \\
dw^{2} &=&w^{1}\wedge w^{5}+w^{2}\wedge w^{6}, \notag\\
dw^{3} &=&-w^{4}\wedge w^{5},
\label{eq:structure}
 \text{\ \ } \\
dw^{4} &=&-w^{3}\wedge w^{4}-w^{4}\wedge w^{6}, \notag \\
dw^{5} &=&w^{3}\wedge w^{5}+w^{5}\wedge w^{6}, \notag \\
dw^{6} &=&w^{4}\wedge w^{5}. \notag
\end{eqnarray}

We remark that changing the basis of invariant forms by the formula

\begin{eqnarray}
\overline{w}^{1} &=&a^{\prime }w^{1}+b^{\prime }w^{2},\text{ \ }\overline{w}%
^{2}=c^{\prime }w^{2},\text{ \ }\overline{w}^{3}=w^{3}+\frac{b^{\prime }}{%
a^{\prime }}w^{5}, \notag \\
\overline{w}^{4} &=&\frac{a^{\prime }}{c^{\prime }}w^{4}-\frac{b^{\prime }}{%
c^{\prime }}w^{3}-\frac{b^{\prime 2}}{a^{\prime }c^{\prime }}w^{5}+\frac{%
b^{\prime }}{c^{\prime }}w^{6},\text{ }\overline{w}^{5}=\frac{c^{\prime }}{%
a^{\prime }}w^{5},
 \label{eq:4.2}
  \text{\ \ } \\
\overline{w}^{6} &=&w^{6}-\frac{b^{\prime }}{a^{\prime }}w^{5}, \notag
\end{eqnarray}

where $a^{\prime },b^{\prime }$ and $c^{\prime }$ denote arbitrary constants
with $a^{\prime }c^{\prime }\neq 0,$ the structure equation \eqref{eq:structure} does not
alter. $\ {\Bigl (}A(2,%
\mathbb{R}
),C=(w^{1},\ldots ,w^{6}){\Bigr )}$ determines a Cartan system.

\subsection{Classification under action of A(2,$%
\mathbb{R}
$)}

The systems to be considered are given by $\Sigma =\langle dx_{2}-f$ $%
(x_{1},x_{2})dx_{1}\rangle $ where $f$ denotes a differentiable function. In
this paper I given some ideas for the classification. Since

\begin{equation*}
dx_{2}-f(x_{1},x_{2})dx_{1}=(x_{5}-x_{3}f(x_{1},x_{2}))w^{1}+(x_{6}-x_{4}f(x_{1},x_{2}))w^{2},
\end{equation*}

and

\begin{equation*}
(x_{6}-x_{4}f(x_{1},x_{2}))^{-1}(dx_{2}-f(x_{1},x_{2})dx_{1}=(x_{6}-x_{4}f(x_{1,}x_{2}))^{-1}(x_{5}-x_{3}f(x_{1},x_{2}))w^{1}+w^{2},
\end{equation*}

then $(x_{6}-x_{4}f)^{-1}(x_{5}-x_{3}f)$ forms a characteristic invariant
system. By using this invariant, we reduce the Cartan system to be
submanifold $M_{0}$ defined by the equation $x_{5}-x_{3}f(x_{1},x_{2})=0.$
The equation to be integrated is now given by $w^{2}=0.$ On the submanilfold 
$M_{0}$ , we have

\begin{eqnarray}
w^{5} &=&aw^{1}+bw^{2}, \notag \\
da &=&2aw^{3}-a\overline{w}^{6}+u_{1}w^{1}+u_{2}w^{2},
 \label{eq:4.3}
  \text{\ \ } \\
db &=&bw^{3}+aw^{4}+(u_{2}-b^{2})w^{1}+u_{3}w^{2}. \notag
\end{eqnarray}

and

\begin{eqnarray}
dw^{1} &=&w^{1}\wedge w^{3}+w^{2}\wedge w^{4}, \notag \\
dw^{2} &=&w^{2}\wedge \overline{w}^{6}, \notag \\
dw^{3} &=&aw^{1}\wedge w^{4}+bw^{2}\wedge w^{4},
 \label{eq:4.4}
  \text{\ \ } \\
dw^{4} &=&-w^{3}\wedge w^{4}-w^{4}\wedge \overline{w}^{6}+bw^{1}\wedge w^{4},
\notag\\
d\overline{w}^{6} &=&-2bw^{2}\wedge w^{4}+u_{3}w^{1}\wedge w^{2}, \notag
\end{eqnarray}

where we put $\overline{w}^{6}=w^{6}-bw^{1}$ and $a,b,u_{1},u_{2},u_{3}$
denote definite functions on $M_{0}.$ These functions are all invariants of
induced Cartan system.

\paragraph{4.2.1 From now on, we shall determine all the equations wich
admit at least 2-dimensional Lie subgroup of $A(2,%
\mathbb{R}
)$ as an invariant group.}

Therefore we suppose always that the forms $w^{1},$ $w^{2}$ are linearly
independent.

\textbf{I. The case }$a=0$. \ From \eqref{eq:4.3} , we have $u_{1}=u_{2}=0$.
Moreover
\begin{equation}\label{eq:4.5}
du_{3}=u_{3}\overline{w}^{6}+u_{3}w^{3}-2bu_{3}w^{1}+cw^{2}
\end{equation}

where $c$ denotes a definite function on $M$.

\textbf{1. }$b=0.$ The manifold $M_{0}$ is given by the maximal integral
manifold of $w^{5}=0.$ Hence we obtain the first type:

\begin{eqnarray*}
\Sigma &:&w^{5}=0, \\
dw^{1} &=&w^{1}\wedge w^{3}+w^{2}\wedge w^{4}, \\
dw^{2} &=&w^{2}\wedge w^{6}, \\
dw^{3} &=&0, \\
dw^{4} &=&-w^{3}\wedge w^{4}-w^{4}\wedge w^{6}, \\
dw^{6} &=&0.
\end{eqnarray*}

Integrating the system $\Sigma ,$ we have the result:

\begin{theorem}
The equation $y^{\prime }=c$ (constant) admits a 5-dimensional Lie subgroup
in $A(2,%
\mathbb{R}
)$ and can be transformed to the equation $y^{\prime }=0$ by an element of $%
A(2,%
\mathbb{R}
).$
\end{theorem}

\textbf{2. }$b\neq 0.$ We can reduce $M_{0}$ to the submanifold $M_{1}$
defined by the equation $b=const.(\neq 0).$ Taking $a^{\prime }=b,$ $%
b^{\prime }=0,$ $c^{\prime }=b$ in \eqref{eq:4.2} we can assume that the constant
is equal to 1: $M_{1}=\left\{ g\in M_{0};\text{ }b(g)=1\right\} .$ From \eqref{eq:4.3},
 \eqref{eq:4.4} and \eqref{eq:4.5} we have on $M_{1}$

\begin{eqnarray}
w^{3} &=&w^{1}-u_{3}w^{2},\text{ \ }w^{5}=w^{2}, \notag \\
du_{3} &=&u_{3}\overline{w}^{6}-u_{3}w^{1}+(c-u_{3}^{2})w^{2}, \notag \\
dw^{1} &=&-u_{3}w^{1}\wedge w^{2}+w^{2}\wedge w^{4},
   \label{eq:4.6}
    \text{\ \ } \\
dw^{2} &=&w^{2}\wedge \overline{w}^{6}, \notag \\
dw^{4} &=&u_{3}w^{2}\wedge w^{4}-w^{4}\wedge \overline{w}^{6}, \notag \\
d\overline{w}^{6} &=&-2w^{2}\wedge w^{4}+u_{3}w^{1}\wedge w^{2}. \notag
\end{eqnarray}

\textbf{2.1. u}$_{3}=0.$ From the second equation of \eqref{eq:4.6} we have $c=0.$
We have thus obtained the second type:

\begin{eqnarray*}
\Sigma &:&w^{3}=w^{1},\text{ \ }w^{5}=w^{2}. \\
dw^{1} &=&w^{2}\wedge w^{4}, \\
dw^{2} &=&w^{1}\wedge w^{2}+w^{2}\wedge w^{6}, \\
dw^{4} &=&-w^{1}\wedge w^{4}-w^{4}\wedge w^{6}, \\
dw^{6} &=&-w^{2}\wedge w^{4}.
\end{eqnarray*}

Integrating the system $\Sigma ,$ we obtain the result:

\begin{theorem}
The equation $y^{\prime }=(x+a)^{-1}(y+b),$ $a,b$ constants, admits a
4-dimensional Lie subgroup in $A(2,%
\mathbb{R}
)$ and can be transformed to $y^{\prime }=x^{-1}y$ by an element of $A(2,%
\mathbb{R}
).$
\end{theorem}

\textbf{2.2. u}$_{3}\neq 0.$ We can reduce $M_{1\text{ }}$to be submanifold $%
M_{2}$ defined by the equation $u_{3}=$ constant $(\neq 0).$ Taking $%
a^{\prime }=1,$ $b^{\prime }=0,$ $c^{\prime }=u_{3}$ in \eqref{eq:4.2}, we can
assume that the constant is equal to 1. From \eqref{eq:4.6} , we have on $M_{2}$

\begin{eqnarray}
w^{3} &=&w^{1}-w^{2},\text{ }w^{5}=w^{2},\text{ }\overline{w}%
^{6}=w^{1}-ew^{2}\text{ \ }(e=c-1), \notag \\
de &=&-3w^{4}+(e-2)w^{1}+rw^{2}, \notag \\
dw^{1} &=&-w^{1}\wedge w^{2}+w^{2}\wedge w^{4},
\label{eq:4.7}
 \text{ \ \ } \\
  dw^{2} &=&-w^{1}\wedge w^{2},\notag \\
  dw^{4} &=&(1-e)w^{2}\wedge w^{4}+w^{1}\wedge w^{4},\notag
\end{eqnarray}

where $r$ denotes a definite function on $M_{2}.$ The equations \eqref{eq:4.7} does
not determines a 3-dimensional Lie group. By using the invariant $e,$ we can
reduce $M_{2.\text{ }}$ Taking $a^{\prime }=1,$ $b^{\prime }=1,$ $c^{\prime
}=\frac{(2-e)}{3}$ in \eqref{eq:4.2}, we can reduce $M_{2}$ to the manifold $%
M_{3}=\left\{ g\in M_{2\text{ }};e(g)=2\right\} .$ From \eqref{eq:4.7} we obtain
on $M_{3}$

\begin{eqnarray}
w^{5} &=&w^{2},\text{ }w^{3}=w^{1}-w^{2},\text{ }\overline{w}%
^{6}=w^{1}-2w^{2},\text{ }w^{4}=\frac{1}{3}rw^{2}, \notag \\
dr &=&2rw^{1}+r_{0}w^{2}
\label{eq:4.8}
\end{eqnarray}
and
\begin{eqnarray}
dw^{1} &=&-w^{1}\wedge w^{2}, \notag \\
dw^{2} &=&-w^{1}\wedge w^{2}.\text{ \ \ \ \ \ \ \ \ \ \ \ \ \ \ }%
\label{eq:4.9}
\end{eqnarray}

Although the equations \eqref{eq:4.9} do not contain any functions, the function $r$
is an invariant of the group. Therefore \eqref{eq:4.8} and \eqref{eq:4.9} determine a
2-dimensional Lie group if and only if $r$ is a constant on $M_{3}.$ In the
case, we have $r=0$ and

\begin{eqnarray*}
\Sigma &:&w^{3}=w^{1}-w^{2},\text{ }w^{4}=0,\text{ }w^{5}=w^{2},\text{ }%
w^{6}=2w^{1}-2w^{2}. \\
dw^{1} &=&-w^{1}\wedge w^{2},\text{ \ }dw^{2}=-w^{1}\wedge w^{2}.
\end{eqnarray*}

Integrating the system $\Sigma ,$ we obtain the result:

\begin{theorem}
Let $f$ be a function satysfying the Clairaut equation

\begin{equation*}
y=xf+\frac{af^{2}+bf+c}{\alpha f+\beta }\text{ \ with }a,\text{ }b,\text{ }c,%
\text{ }\alpha ,\text{ }\beta \text{ constants.}
\end{equation*}
\end{theorem}

Then $y^{\prime }=f(x,y)$ admits a 2-dimensional Lie group in $A(2,%
\mathbb{R}
)$ and can be transformed to $y^{\prime }=x-(x^{2}-2y)^{\frac{1}{2}}$ by an
element $A(2,%
\mathbb{R}
).$

\textbf{II. The case }$a\neq 0.$ We go back to the manifold $M_{0}.$ Suppose
that $b\neq 0$ on $M_{0}$ . Consider the submanifold $N_{0}$ defined by the
equation $b=const$ $(\neq 0).$ By setting 
\begin{eqnarray*}
\overline{w}^{1} &=&bw^{1}+\frac{b^{2}}{a}w^{2},\text{ \ }\overline{w}^{2}=%
\frac{b}{a}w^{2},\text{ \ }\overline{w}^{3}=w^{3}+\frac{b}{a}w^{5}, \\
\overline{w}^{4} &=&aw^{4}-bw^{3}-\frac{b^{2}}{a}w^{5}+bw^{6},\text{ }%
\overline{w}^{5}=\frac{1}{a}w^{5},\text{ \ }\overline{w}^{6}=w^{6}-\frac{b}{a%
}w^{5},
\end{eqnarray*}

we can assume that $a=1,b=0$ on $M_{0}.$ Therefore we have only to examine
the case $a\neq 0,$ $b=0$ on $M_{0}.$ We can reduce $M_{0}$ to the
submanifold

\begin{equation*}
N_{0}=\left\{ g\in M_{0};\text{ }a(g)=1,b(g)=0\right\}
\end{equation*}

on which we have

\begin{eqnarray}
w^{4} &=&-u_{2}w^{1}-u_{3}w^{2},\text{ }w^{5}=w^{1}, \text{ }%
w^{6}=2w^{3}+u_{1}w^{1}+u_{2}w^{2},\notag \\
dw^{1} &=&w^{1}\wedge w^{3}+u_{2}w^{1}\wedge w^{2},
 \label{eq:4.10}
  \text{ \ \ } \\
dw^{2} &=&2w^{2}\wedge w^{3}-u_{1}w^{1}\wedge w^{2}, \notag \\
dw^{3} &=&-u_{3}w^{1}\wedge w^{2}. \notag 
\end{eqnarray}

By differentiating \eqref{eq:4.9} we obtain

\begin{eqnarray}
du_{1} &=&u_{1}w^{3}+v_{1}w^{1}+v_{2}w^{2}, \notag \\
du_{2} &=&2u_{2}w^{3}+v_{3}w^{1}+v_{4}w^{2},
\label{eq:4.11}
 \text{ \ \ } \\
du_{3} &=&3u_{3}w^{3}+v_{5}w^{1}+v_{6}w^{2}, \notag \\
0 &=&-v_{5}+v_{4}+2(u_{1}u_{3}-u_{2}^{2}), \notag \\
0 &=&v_{2}-v_{3}+3u_{3}. \notag
\end{eqnarray}

\textbf{1. }$u_{1}=u_{2}=u_{3}=0.$ From \eqref{eq:4.10}, we have

\begin{eqnarray*}
\Sigma &:&\text{ }w^{4}=0,\text{ }w^{5}=w^{1},\text{ }w^{6}=2w^{3}. \\
dw^{1} &=&w^{1}\wedge w^{3}, \\
dw^{2} &=&2w^{2}\wedge w^{3}, \\
dw^{3} &=&0.
\end{eqnarray*}

Integrating the system $\Sigma ,$ we obtain the result:

\begin{theorem}
The equations of this type can be transformed to $y^{\prime }=x$ by an
element of $A(2,%
\mathbb{R}
)$ and admit a 3-dimensional Lie group in $A(2,%
\mathbb{R}
).$
\end{theorem}

\textbf{2. }$u_{1}\neq 0,$ $u_{2}=0,$ $u_{3}=0.$ We reduce $N_{0}$ to the
submanifold defined by the equation $u_{1}=const.(\neq 0).$ From \eqref{eq:4.11}
we have $v_{2}=v_{3}=v_{4}=v_{5}=v_{6}=0$ and hence

\begin{eqnarray*}
\Sigma &:&\text{ }w^{3}=-\frac{v_{1}}{u_{1}}w^{1},\text{ }w^{4}=0,\text{ }%
w^{5}=w^{1},\text{ }w^{6}=(u_{1}-\frac{2v_{1}}{u_{1}})w^{1}. \\
dw^{1} &=&0, \\
dw^{2} &=&-(u_{1}-\frac{2v_{1}}{u_{1}})w^{1}\wedge w^{2}, \\
dv_{1} &=&0(\func{mod}w^{1}).
\end{eqnarray*}

If these equations determine a 2-dimensional Lie group, $v_{1}$ must be a
constant. In this case, integrating the system $\Sigma ,$ we obtain the
result:

\begin{theorem}
All the equations in this case are transformed by an element of $A(2,%
\mathbb{R}
)$ to one the following three types

i) $y^{\prime }=\log x,$ \ ii) $y^{\prime }=e^{x},$ \ iii) $y^{\prime
}=x^{a} $ \ $(a$ $const.\neq 0,1),$ which admit 2-dimensional Lie group in $%
A(2,%
\mathbb{R}
).$
\end{theorem}

\textbf{3.\ }$u_{1}=0,$ $u_{2}\neq 0,$ $u_{3}=0.$ We reduce $N_{0}$ to the
submanifold defined by the equation $u_{2}=const.(\neq 0).$ From \eqref{eq:4.11},
we have $v_{4}=2u_{2}^{2}$ and otherwise $v_{i}=0$ and hence

\begin{eqnarray*}
\Sigma &:&w^{3}=-u_{2}w^{2},\text{ }w^{4}=-u_{2}w^{1},\text{ }w^{5}=w^{1},%
\text{ }w^{6}=-u_{2}w^{2}. \\
dw^{1} &=&0,\text{ }dw^{2}=0.
\end{eqnarray*}

\textbf{3.1. The case }$u_{2}$ $\gvertneqq $ $0$. \ We can assume always $%
u_{2}=1.$ Integrating the system $\Sigma ,$ we obtain the result:

\begin{proposition}
All the equations in this type are transformed by an element of $A(2,%
\mathbb{R}
)$ to the equation $y^{\prime }=-xy^{-1},$ which admits a 2-dimensional Lie
group in $A(2,%
\mathbb{R}
).$
\end{proposition}

\textbf{3.2. The case }$u_{2}\lvertneqq 0.$ We can assume $u_{2}=-1.$
Integrating the system $\Sigma $ we obtain the result:

\begin{proposition}
All equations in this case are transformed by an element of $A(2,%
\mathbb{R}
)$ to the equation $y^{\prime }=-x^{-1}y,$ which admits a 2-dimensional Lie
group in $A(2,%
\mathbb{R}
).$
\end{proposition}

\textbf{4. }$u_{1}=u_{2}=0,$ $u_{3}\neq 0.$ From \eqref{eq:4.11} we have $v_{i}=0,$
$i=1,2,3,4.$ Since $v_{2}-v_{3}+3u_{3}=0,$ this contradicts the assumption $%
u_{3}\neq 0.$

\textbf{5. }$u_{1}=0,$ $u_{2}\neq 0,$ $u_{3}\neq 0.$ We reduce $N_{0}$ to
the submanifold defined by the equations $u_{2}$ $=const.(\neq 0),$ $%
u_{3}=const.(\neq 0).$ We can assume $u_{3}=4.$ From \eqref{eq:4.11}, we have $%
v_{1}=v_{2}=0,$ $v_{3}=12,$ $v_{4}-v_{5}=2u_{2}^{2}$ and

\begin{eqnarray}
0 &=&2u_{2}w^{3}+12w^{1}+v_{4}w^{2},
 \label{eq:4.12}
  \text{\ \ } \\
0 &=&12w^{3}+v_{5}w^{1}+v_{6}w^{2}. \notag
\end{eqnarray}

By this equations, we obtain

\begin{equation*}
v_{4}=\frac{72}{u_{2}}+2u_{2}^{2},\text{ \ }v_{5}=\frac{72}{u_{2}},\text{ \ }%
v_{6}=\frac{6}{u_{2}}(\frac{72}{u_{2}}+2u_{2}^{2}).
\end{equation*}

Substituting these values to \eqref{eq:4.12} we have

\begin{equation*}
w^{3}=-\frac{6}{u_{2}}w^{1}-(\frac{36}{u_{2}^{2}}+u_{2})w^{2}.
\end{equation*}

Substituting this equation to the last equation of \eqref{eq:4.10} we have $%
u_{2}=-3.$ Hence we have

\begin{eqnarray*}
\Sigma &:&\text{ }w^{3}=2w^{1}-w^{2},\text{ \ }w^{4}=3w^{1}-4w^{2},\text{ \ }%
w^{5}=w^{1},\text{ \ }w^{6}=2w^{3}-3w^{2}. \\
dw^{1} &=&-4w^{1}\wedge w^{2},\text{ \ }dw^{2}=-4w^{1}\wedge w^{2}.
\end{eqnarray*}

Integrating the system $\Sigma ,$ we obtain the result:

\begin{theorem}
The equation of this type can be transformed by an element of $A(2,%
\mathbb{R}
)$ to $y^{\prime }=x^{-1}y+x^{4},$ which admits a 2-dimensional Lie group in 
$A(2,%
\mathbb{R}
).$
\end{theorem}

\textbf{6. }$u_{1}\neq 0,$ $u_{2}=0,$ $u_{3}\neq 0.$ \ We reduce $N_{0}$ to
the submanifold defined by the equations $u_{1}=const.(\neq 0),$ $%
u_{3}=const.(\neq 0).$ We can assume $u_{1}=-3.$ By the same argument as in
the case \textbf{5}, we have

\begin{eqnarray*}
\Sigma &:&\text{ }w^{3}=2w^{1}-u_{3}w^{2},\text{ }w^{4}=-u_{3}w^{2},\text{ }%
w^{5}=w^{1},\text{ }w^{6}=w^{1}-2u_{3}w^{2}. \\
dw^{1} &=&-u_{3}w^{1}\wedge w^{2},\text{ \ }dw^{2}=-w^{1}\wedge w^{2}.
\end{eqnarray*}

Integrating the system $\Sigma ,$ we obtain the result:

\begin{theorem}
i) If $u_{3}=1,$ the equations of this type are transformed by an element of 
$A(2,%
\mathbb{R}
)$ to $y^{\prime }=x^{-1}y+x.$

ii) If $u_{3}\neq 1,$ the equations of this type are transformed by an
element of $A(2,%
\mathbb{R}
)$ to the equation

\begin{equation*}
y^{\prime }=\frac{-u_{3}x+\sqrt[2]{2(1-u_{3})y+u_{3}x^{2}}}{1-u_{3}}.
\end{equation*}
\end{theorem}

\textbf{7. }$u_{1}\neq 0,$ $u_{2}\neq 0,$ $u_{3}\neq 0.$ By the same
argument as in the case \textbf{5, }we have

\begin{eqnarray*}
\Sigma &:&\text{ }w^{3}=-\frac{1}{2}u_{1}w^{1}-u_{2}w^{2},\text{ }%
w^{4}=-u_{2}w^{1},\text{ }w^{5}=w^{1},\text{ }%
w^{6}=2w^{3}+u_{1}w^{1}+u_{2}w^{2}. \\
dw^{1} &=&dw^{2}=0.
\end{eqnarray*}

Integrating the system $\Sigma ,$ we have the result:

\begin{proposition}
The equations of this type are transformed by an element of $A(2,%
\mathbb{R}
)$ to the equation

\begin{equation*}
y^{\prime }=\frac{-2x+u_{1}y}{2u_{2}y}\text{ \ .}
\end{equation*}
\end{proposition}

\textbf{8. }$u_{1}\neq 0,$ $u_{2}\neq 0,$ $u_{3}\neq 0.$ \ By the same
argument as in the case \textbf{5 , }we can deternime $v_{i}$ $,1\leq i\leq
6.$ In particular, we have

\begin{equation*}
v_{1}=\frac{u_{1}(6u_{2}-2u_{1}^{2}+u_{1}u_{2}^{2})}{2u_{2}^{2}-3u_{1}},%
\text{ \ }v_{2}=\frac{u_{1}(18-4u_{1}u_{2}+2u_{2}^{3})}{2u_{2}^{2}-3u_{1}}.
\end{equation*}

If $2u_{2}^{2}=3u_{1},$ then $u_{1}=6,$ $u_{2}=3,$ $u_{3}=2,$ $v_{1}=2$ and $%
v_{2}=1.$ We have

\begin{equation*}
w^{3}=-\frac{v_{1}}{u_{1}}w^{1}-\frac{v_{2}}{u_{1}}w^{2}.
\end{equation*}

Substituting this relation to the last equation of \eqref{eq:4.10} we have a
certain algebraic equation with respect to the quantities $u_{1},u_{2}.$
Hence we obtain

\begin{eqnarray*}
\Sigma &:&\text{ }w^{3}=-\frac{v_{1}}{u_{1}}w^{1}-\frac{v_{2}}{u_{1}}w^{2},%
\text{ }w^{4}=-u_{2}w^{1}-u_{3}w^{2},\text{ }w^{5}=w^{1},\text{ }%
w^{6}=2w^{3}+u_{1}w^{1}+u_{2}w^{2}. \\
dw^{1} &=&(u_{2}-\frac{v_{2}}{u_{1}})w^{1}\wedge w^{2},\text{ \ \ }%
dw^{2}=-(u_{1}-\frac{2v_{1}}{u_{1}})w^{1}\wedge w^{2}.
\end{eqnarray*}

Intregating the system $\Sigma ,$ we obtain the result:

\begin{theorem}
The equations of this type are transformed by an element of $A(2,%
\mathbb{R}
)$ to $y^{\prime }=x^{-1}y+x^{a}$ $(a$ $const\neq 4,1,0,-1),$ which admits a
2-dimensional Lie group in $A(2,%
\mathbb{R}
).$
\end{theorem}

\paragraph{4.2.2. \ As for the determination of the equations admitting a
1-dimensional Lie group in $A(2,%
\mathbb{R}
),$ we can use the method developed in section 3.}

Here is the table of the standard forms and the invariant groups. We denote
by $a$ the parameter of a 1-dimensional Lie group.
\begin{center}
\begin{tabular}{|l|l|}
 \hline\\
Standard Forms & Invariant Groups\\
\hline
\hline\\
$y'=F(x)$ & $X=x$, $Y=y+a$\\
\hline\\
$y'=\frac{y}{x}F(\frac{y^{r}}{x^{s}})$ & $X=a^{r}x$, $Y=a^{s}y$\\
\hline\\
$y'=yF(ye^{-x})$ & $X=x+a$, $Y=e^{a}y$\\
\hline\\
$y'=\frac{y}{x}$\ $+F(x)$ & $X=x, Y=ax+y$\\
\hline\\
$y'=\frac{y}{x}+F(xe^{r\frac{y}{x}})$ & $X=e^{ar}x$, $Y=ae^{ar}x+e^{ar}y$\\
\hline\\
$y'=\frac{y-xF(x^{2}+y^{2})}{x+yF(x^{2}+y^{2})}$ & \text{rotation group}\\
\hline\\
$\frac{y-xy'}{x+yy'}=F\left(\frac{y-x\tan(r\log\sqrt[2]{x^{2}+y^{2}})}{x+y\tan(r\log\sqrt[2]{x^{2}+y^{2}})}\right)$ & \text{1-dimensional conformal transformation group}\\
\hline

\end{tabular} 
\end{center}

\end{document}